\documentclass[12pt,reqno]{amsart}

\theoremstyle{plain}    
\newtheorem{thm}{Theorem}[section]
\numberwithin{figure}{section} 
\theoremstyle{plain}    
\newtheorem{prop}[thm]{Proposition}

\newtheorem{defi}[thm]{Definition}

\newtheorem{examples}[thm]{Examples}

\newtheorem{remark}[thm]{Remark}

\usepackage{amscd,amssymb,comment,euscript}
\newcommand\dif{\mbox{\it d}}

\newcommand\DT{\operatorname{DT}}

\newcommand\eps{\epsilon}

\newcommand\EEu{{\EuScript E}}                   

\newcommand\Ft{{\widetilde F}}

\newcommand\HEu{{\EuScript H}}                   

\newcommand\Nats{{\mathbf N}}

\newcommand\Thetat{{\widetilde\Theta}}

\newcommand\ttild{{\tilde t}}

\begin{document}

\pagestyle{myheadings}

\title{Generating functions for moments of the quasi--nilpotent DT--operator}
 
\author{Ken Dykema}

\author{Catherine Yan}

\address{\hskip-\parindent
Department of Mathematics\\
Texas A\&M University\\
College Station TX 77843--3368, USA}

\email{kdykema@math.tamu.edu}
\email{cyan@math.tamu.edu}

\thanks{K.D.\ supported in part by NSF grant DMS--0070558.
He thanks also the Mathematical Sciences Research Institute,
where part of this work was carried out.
Research at MSRI is supported in part by NSF grant DMS--9701755.
C.Y.\ supported in part by NSF grant DMS-0070574.
She is also
partly supported by NSF grant DMS-9729992 through the Institute
for Advanced Study.}

\begin{abstract}
We prove a recursion formula for generating functions of certain renormalizations of
$*$--moments of the $\DT(\delta_0,1)$--operator $T$,
involving an operation $\odot$ on formal power series and a transformation $\EEu$
that converts $\odot$ to usual multiplication.
This recursion formula is used to prove that all of these generating functions
are rational functions, and to find a few of them explicitly.
\end{abstract}

\date{Jan. 10, 2002}

\maketitle

\markboth{\tiny T-moments}{\tiny T-moments}
 
\section*{Introduction}

In combinatorics, one of the most useful methods for studying a sequence  is to 
give its generating functions. The two most common types of generating functions are
ordinary generating functions $\sum f(n)x^n$, and exponential generating functions
$\sum f(n)x^n/n!$. In this paper we study the $*$--moment generating functions --- a
family of multivariable power series $F_n$--- of a particular operator $T$ that arose in 
the theory of free probability. We prove that $F_n$'s are all rational 
by applying a linear transformation between these two types of generating 
functions. 

The central object of this paper is the collection of $*$--moments of a particular bounded
operator $T$ on Hilbert space, which was
constructed in~\cite{DH} and which is a candidate for an operator
without a nontrivial hyperinvariant subspace.
(A hyperinvariant subspace of an operator $T$ on a Hilbert space $\HEu$ is a closed subspace
$\HEu_0\subseteq\HEu$ that is invariant under every operator $S$ that commutes with $T$,
i.e. $S(\HEu_0)\subseteq\HEu_0$.
It is an open problem
whether every operator on Hilbert space that is not a multiple
of the identity has a nontrivial hyperinvariant subspace.)
The von Neumann algebra generated by $T$
has a unique normalized trace $\tau$,
and by the $*$--moments of $T$ we mean the values
\[
M(k_1,\ell_1,\dots,k_n,\ell_n)=\tau\big((T^*)^{k_1}T^{\ell_1}\dots(T^*)^{k_n}T^{\ell_n}\big),
\]
with $n\in\Nats$, $k_1,\dots,k_n,\ell_1,\dots,\ell_n\in\Nats\cup\{0\}$.
These $*$--moments determine a representation of $T$ on a Hilbert space, (which can be shown to be bounded, see~\cite{DH}),
via the construction
of Gelfand, Naimark and Segal, (cf.~\cite{KR})
and hence they encode all essential properties of the operator.
Our effort to understand the $*$--moments is part of an attempt better to
understand the operator $T$.

To be precise, $T$ is a $\DT(\delta_0,1)$--operator, in the class of DT--operators
constructed in~\cite{DH}.
$T$ can be realized as the limit in $*$-moments of strictly upper triangular random
matrices having i.i.d.\ complex gaussian entries.
In~\cite{DH}, it was proved that $T$ is quasi--nilpotent, i.e.\ has spectrum $\{0\}$.
It was shown that $M(k_1,\ell_1,\dots,k_n,\ell_n)=0$
if $k_1+\dots+k_n\ne\ell_1+\dots+\ell_n$, and also a recursion formula for $M$
was proved (cf.~\cite[Theorem 9.5]{DH}).
It was conjectured in~\cite{DH} that
\begin{equation}
\label{eq:conj}
\tau\big(((T^*)^kT^k)^n)=\frac{n^{nk}}{(nk+1)!}\qquad(k,n\in\Nats).
\end{equation}

In this paper, we will consider the quantities
\[
N(k_1,\ell_1,\dots,k_n,\ell_n)=\begin{cases}
0&\text{if }
\begin{aligned}[t]
&k_1+\dots+k_n \\ &\ne\ell_1+\dots+\ell_n
\end{aligned} \\
(m+1)!\,M(k_1,\ell_1,\dots,k_n,\ell_n)&\text{if }
\begin{aligned}[t]
m:=&k_1+\dots+k_n \\
=&\ell_1+\dots+\ell_n
\end{aligned}
\end{cases}
\]
(which must be nonnegative integers by~\cite[Lemma 2.4]{DH}),
and their generating functions
\begin{equation}
\label{eq:Fdef}
F_n(z_1,w_1,\dots,z_n,w_n)=\sum_{\substack{k_1,\dots,k_n\ge0 \\ \ell_1,\dots,\ell_n\ge0}}
N(k_1,\ell_1,\dots,k_n,\ell_n)z_1^{k_1}w_1^{\ell_1}\dots z_n^{k_n}w_n^{\ell_n}.
\end{equation}
We will prove a recursion formula for the functions $F_n$, involving an operation $\odot$
on power series, which can be decsribed as multiplication with reweighting of homogeneous parts.
In order to compute these generating functions,
we will define a linear transformation of formal
power series, called the $\EEu$--transform,
which maps these power series to exponential generating functions of a new variable $q$, and
converts the operation $\odot$  to the usual multiplication.
This will allow us to prove that $F_n$ is a rational function for all $n$,
and to give an algorithm for computing $F_n$.

In~\S\ref{sec:fps}, the operation $\odot$ and the transform $\EEu$ on formal power
series are introduced, and several properties are proved.
In~\S\ref{sec:genfun}, recursion formulas for $N(k_1,\ell_1,\ldots,k_n,\ell_n)$
and for the generating functions $F_n$ are proved.
The latter are shown to be rational and several examples are computed.
In~\S\ref{sec:ci}, we use the method of contour integration to verify the conjecture~\eqref{eq:conj}
in the case $n=2$.

\section{Formal power series}
\label{sec:fps}

In this section, fix variables $x_1,x_2,\dots$ and
fix $N\in\Nats$.
(In the next section we will use only the case $N=2$.)
Let $\Theta=\Theta_N$
be the algebra of all formal power series in $x_1,x_2,\dots$ of the form
\begin{equation}
\label{eq:fps}
f(x_1,\dots,x_n)=\sum_{k_1,\dots,k_n\ge0}c_{k_1,\dots,k_n}x_1^{k_1}\dots x_n^{k_n}
\end{equation}
for some $n\in\Nats$, and
such that $c_{k_1,\dots,k_n}=0$ whenever $k_1+\dots+k_n$ is not divisible by $N$.
Given a formal power series as in~\eqref{eq:fps} and given $k\in\Nats\cup\{0\}$,
write $f^{(k)}$ for its $Nk$--homogeneous part:
\begin{equation}
\label{eq:fhomog}
f^{(k)}(x_1,\dots,x_n)=\sum_{k_1+\dots+k_n=Nk}
c_{k_1,\dots,k_n}x_1^{k_1}\dots x_n^{k_n}.
\end{equation}

\begin{defi}\rm
\label{def:dtimes}
Let $\odot=\odot_N$ be the binary operation on $\Theta$ given by
\[
f\odot g=\sum_{k,\ell\ge0}\binom{k+\ell}k
f^{(k)}g^{(\ell)}.
\]
\end{defi}

Note that $\odot$ is bilinear, commutative and associative.
If $f_i\in\Theta$ ($i\in\{1,\dots,p\}$) then
\[
(f_1\odot f_2\odot\dots\odot f_p)
=\sum_{k_1,\dots,k_p\ge0}\binom{k_1+\dots+k_p}{k_1,k_2,\dots,k_p}
\prod_{i=1}^pf_i^{(k_i)}.
\]

Let $\Thetat$ be the set of all formal power series in variables $x_1,x_2,\dots$ and
the additional variable~$q$.

\begin{defi}\rm
\label{def:Etrans}
The $\EEu$--transform $\EEu=\EEu_N:\Theta\to\Thetat$ is the map
given by
\[
\EEu f=\sum_{k\ge0}\frac{q^k}{k!}f^{(k)}.
\]
\end{defi}

Note that $\EEu$ is linear.

\begin{prop}
\label{prop:basic}
\renewcommand{\theenumi}{\roman{enumi}}
\renewcommand{\labelenumi}{(\theenumi)}
\begin{enumerate}

\item
\label{enum:fdtg}
If $f,g\in\Theta$ then
$\EEu(f\odot g)=(\EEu f)(\EEu g)$.

\item
\label{enum:af}
If $f\in\Theta$ and if $a\in\Theta$ is a homogeneous polynomial
of degree~$N$ then
$\frac{\dif}{\dif q}\big(\EEu(af)\big)=a\,\EEu f$,
i.e.
\[
\EEu(af)(x_1,x_2,\dots,q)=a\int_0^q\EEu f(x_1,x_2\dots,t)\dif t.
\]

\item
\label{enum:1u}
If $r\in\Nats$ and if $u_1,\dots,u_r\in\Theta$
are distinct homogeneous polynomials of degree $N$, then
\[
\EEu\bigg(\frac1{(1-u_1)(1-u_2)\dots(1-u_r)}\bigg)=
\sum_{i=1}^r\frac{u_i^{r-1}e^{qu_i}}{\prod_{j\ne i}(u_i-u_j)},
\]
which in the case $r=1$ means
$\EEu(1/(1-u_1))=e^{qu_1}$.

\item
\label{enum:a1u}
Let $r\in\Nats$ and $k\in\{0,1,\dots,r-1\}$.
If $u_1,\dots,u_r,a_1,\dots,a_k\in\Theta$
are homogeneous polynomials of degree $N$ with
$u_1,\dots,u_r$ distinct,
then
\[
\EEu\bigg(\frac{a_1\dots a_k}{(1-u_1)(1-u_2)\dots(1-u_r)}\bigg)=
a_1\dots a_k\sum_{i=1}^r\frac{u_i^{r-k-1}e^{qu_i}}{\prod_{j\ne i}(u_i-u_j).}
\]

\end{enumerate}
\end{prop}
\begin{proof}
For~(\ref{enum:fdtg}),
\[
\EEu(f\odot g)=\sum_{r\ge0}\frac{q^r}{r!}
\sum_{k+\ell=r}\binom rkf^{(k)}g^{(\ell)}
=\sum_{k,\ell\ge0}\bigg(\frac{q^k}{k!}f^{(k)}\bigg)
\bigg(\frac{q^\ell}{\ell!}g^{(\ell)}\bigg)
=(\EEu f)(\EEu g).
\]

For~(\ref{enum:af}),
\[
\frac\dif{\dif q}\EEu(af)=\frac\dif{\dif q}\EEu\bigg(\sum_{k=0}^\infty af^{(k)}\bigg)
=a\frac\dif{\dif q}\sum_{k=0}^\infty \frac{q^{k+1}}{(k+1)!}f^{(k)}
=a\sum_{k=0}^\infty \frac{q^k}{k!}f^{(k)}
=a\,\EEu f.
\]

For~(\ref{enum:1u}),
\[
\EEu\bigg(\frac1{(1-u_1)\dots(1-u_r)}\bigg)=
\sum_{k=0}^\infty\frac{q^k}{k!}\sum_{k_1+\dots+k_p=k}u_1^{k_1}\dots u_r^{k_r}.
\]
On the other hand,
\[
\sum_{i=1}^r\frac{u_i^{r-1}e^{qu_i}}{\prod_{j\ne i}(u_i-u_j)}
=\sum_{k=0}^\infty\frac{q^k}{k!}\sum_{i=1}^r\frac{u_i^{k+r-1}}{\prod_{j\ne i}(u_i-u_j).}
\]
Hence it will suffice to show
\begin{equation}
\label{eq:sumutoshow}
\sum_{k_1+\dots+k_r=k}u_1^{k_1}\dots u_r^{k_r}
=\sum_{i=1}^r\frac{u_i^{k+r-1}}{\prod_{j\ne i}(u_i-u_j)}
\end{equation}
for all $k\ge0$.
We will proceed by induction on $r$,
treating $u_1,\dots,u_r$ in~\eqref{eq:sumutoshow} as any $n$ distinct elements in a field.
The case $r=1$ is clear.
If $r=2$, then we have
\begin{align*}
\sum_{k_1+k_2=k}u_1^{k_1}u_2^{k_2}&=
u_2^k\sum_{k_1=0}^k\bigg(\frac{u_1}{u_2}\bigg)^{k_1}
=u_2^k\frac{(u_1/u_2)^{k+1}-1}{(u_1/u_2)-1} \\
&=\frac{u_1^{k+1}-u_2^{k+1}}{u_1-u_2}
=\frac{u_1^{k+1}}{u_1-u_2}+\frac{u_2^{k+1}}{u_2-u_1},
\end{align*}
as required.
Assume $r\ge3$.
Using the induction hypothesis, we have
\begin{align*}
\sum_{k_1+\dots+k_r=k}u_1^{k_1}\dots u_r^{k_r}
&=u_1^k\sum_{k'=0}^k\quad\sum_{k_2+\dots+k_r=k'}\bigg(\frac{u_2}{u_1}\bigg)^{k_2}
 \bigg(\frac{u_3}{u_1}\bigg)^{k_3}\dots \bigg(\frac{u_r}{u_1}\bigg)^{k_r} \\
&=u_1^k\sum_{k'=0}^k\sum_{i=2}^r\frac{(u_i/u_1)^{k'+r-2}}{\prod_{j\ne1,i}\big((u_i/u_1)-(u_j/u_1)\big)} \\
&=u_1^k\sum_{i=2}^r\frac{\big((u_i/u_1)^{k+1}-1\big)(u_i/u_1)^{r-2}}{\big((u_i/u_1)-1\big)
 \prod_{j\ne1,i}\big((u_i/u_1)-(u_j/u_1)\big)} \\
&=\sum_{i=2}^r\frac{u_i^{k+r-1}}{\prod_{j\ne i}(u_i-u_j)}
 -u_1^{k+1}\sum_{i=2}^r\frac{u_i^{r-2}}{\prod_{j\ne i}(u_i-u_j).}
\end{align*}
We must show
\[
-\sum_{i=2}^r\frac{u_i^{r-2}}{\prod_{j\ne i}(u_i-u_j)}
=\frac{u_1^{r-2}}{\prod_{j\ne1}(u_1-u_j)}
\]
or, equivalently,
\begin{equation}
\label{eq:vddtoshow}
\sum_{i=1}^ru_i^{r-2}\prod_{\substack{1\le k<\ell\le r \\ k,\ell\ne i}}(u_\ell-u_k)=0.
\end{equation}
Using the Vandermonde determinant
\[
\det\left(\begin{matrix}
        1 &         1 & \dots &             1 \\
      u_1 &       u_2 & \dots &       u_{r-1} \\
   \vdots &    \vdots & \vdots &        \vdots \\
u_1^{r-2} & u_2^{r-2} & \dots & u_{r-1}^{r-2}
\end{matrix}\right)=\prod_{1\le k<\ell\le r-1}(u_\ell-u_k)
\]
and expanding the determinant
\[
\det\left(\begin{matrix}
u_1^{r-2} & u_2^{r-2} & \dots & u_r^{r-2} \\
        1 &         1 & \dots &         1 \\
      u_1 &       u_2 & \dots &       u_r \\
   \vdots &    \vdots & \vdots &    \vdots \\
u_1^{r-2} & u_2^{r-2} & \dots & u_r^{r-2}
\end{matrix}\right)
\]
along the first row, we find the latter is equal to the LHS of~\eqref{eq:vddtoshow},
but is also clearly equal to zero.

To prove~(\ref{enum:a1u}), use~(\ref{enum:1u}) and iterate~(\ref{enum:af}).
\end{proof}

\begin{remark}\rm
\label{rem:E-1}
If one knows that $h\in\Thetat$ is in the image of $\EEu$, then to find
$\EEu^{-1}(h)$, simply replace every occurrence of $q^k$ by $k!$.
This can sometimes be used to compute $f\odot g$.
\end{remark}

We now turn to some cases of $f\odot g$ which will be used in~\S\ref{sec:genfun}.

\begin{prop}
\label{prop:11v}
If $m\in\Nats$, $k\in\{0,1,\dots,m-1\}$ and if $a_1,\dots,a_k,u_1,\dots,u_m,v\in\Theta$
are homogeneous polynomials of degree $N$ such that $u_1,\dots,u_m$ are distinct, then
\[
\frac{a_1\dots a_k}{(1-u_1)\dots(1-u_m)}\odot\frac1{(1-v)}
=\frac{a_1\dots a_k(1-v)^{m-k-1}}{(1-u_1-v)(1-u_2-v)\dots(1-u_m-v)}.
\]
\end{prop}
\begin{proof}
Using Proposition~\ref{prop:basic}, we have
\begin{align*}
&\frac{a_1\dots a_k}{(1-u_1)\dots(1-u_m)}\odot\frac1{(1-v)} \\
&\qquad\overset\EEu\mapsto\quad
\bigg(a_1\dots a_k\sum_{i=1}^m
\frac{u_i^{m-k-1}e^{qu_i}}{\prod_{j\ne i}(u_i-u_j)}\bigg)e^{qv} \displaybreak[2] \\
&\qquad=\quad a_1\dots a_k\sum_{i=1}^m
\frac{\big((u_i+v)-v\big)^{m-k-1}e^{q(u_i+v)}}{\prod_{j\ne i}(u_i-u_j)} \displaybreak[2] \\
&\qquad=\quad a_1\dots a_k\sum_{p=0}^{m-k-1}(-1)^p\binom{m-k-1}pv^p\sum_{i=1}^m
\frac{(u_i+v)^{m-k-p-1}e^{q(u_i+v)}}{\prod_{j\ne i}\big((u_i+v)-(u_j+v)\big)} \displaybreak[2] \\
&\qquad\overset{\EEu^{-1}}\mapsto
\quad\sum_{p=0}^{m-k-1}(-1)^p\binom{m-k-1}pv^p\frac{a_1\dots a_k}{\prod_i(1-u_i-v)} \displaybreak[2] \\
&\qquad=\quad\frac{a_1\dots a_k(1-v)^{m-k-1}}{(1-u_1-v)(1-u_2-v)\dots(1-u_m-v)}.
\end{align*}
\end{proof}

\begin{prop}
\label{prop:adtb}
Let $m,n\in\Nats$ and let $k\in\{0,1,\dots,m-1\}$, $\ell\in\{0,1,\dots,n-1\}$.
Let $u_1,\dots,u_m,a_1,\dots,a_k,v_1\dots,v_n,b_1,\dots,b_\ell\in\Theta$
be homogeneous polynomials of degree $N$ such that $(u_i+v_j)_{1\le i\le m,\,1\le j\le n}$
is a family of $mn$ distinct polynomials.
Then
\[
\frac{a_1\dots a_k}{(1-u_1)\dots(1-u_m)}\odot\frac{b_1\dots b_\ell}{(1-v_1)\dots(1-v_n)}
=\frac{a_1\dots a_k\,b_1\dots b_\ell\;P^{k,\ell}_{m,n}}{\prod_{i,j}(1-u_i-v_j),}
\]
where $P^{k,\ell}_{m,n}=P^{k,\ell}_{m,n}(u_1,\dots,u_m,v_1,\dots,v_n)$ is a polynomial in the
variables $(u_i)_1^m$ and $(v_j)_1^n$, having degree
bounded above by $mn-k-\ell-1$.
\end{prop}
\begin{proof}
Using Proposition~\ref{prop:basic} and Remark~\ref{rem:E-1}, we have
\begin{align*}
\frac{a_1\dots a_k}{\prod_i(1-u_1)}&\odot\frac{b_1\dots b_\ell}{\prod_j(1-v_j)} \\
\overset\EEu\mapsto&\quad a_1\dots a_k\,b_1\dots b_\ell
 \bigg(\sum_{i=1}^m\frac{u_i^{m-k-1}e^{qu_i}}{\prod_{i'\ne i}(u_i-u_{i'})}\bigg)
 \bigg(\sum_{j=1}^n\frac{v_j^{n-\ell-1}e^{qv_i}}{\prod_{j'\ne j}(v_j-v_{j'})}\bigg) \\
\overset{\EEu^{-1}}\mapsto&\quad a_1\dots a_k\,b_1\dots b_\ell
 \sum_{i,j}\frac{u_i^{m-k-1}v_j^{n-\ell-1}}
 {(1-u_i-v_j)\big(\prod_{i'\ne i}(u_i-u_{i'})\big)\big(\prod_{j'\ne j}(v_j-v_{j'})\big)} \\
=&\quad\bigg(\frac{a_1\dots a_k\,b_1\dots b_\ell}{\prod_{i,j}(1-u_i-v_j)}\bigg)\;
 \frac Q{\big(\prod_{p<q}(u_p-u_q)\big)\big(\prod_{r<s}(v_r-v_s)\big)},
\end{align*}
where $Q=Q(u_1,\dots,u_m,v_1,\dots,v_n)$ is the polynomial
\begin{equation}
\label{eq:Q}
\begin{aligned}
Q=\sum_{i,j}(-1)^{i+j}u_i^{m-k-1}v_j^{n-\ell-1}&
 \bigg(\prod_{(i',j')\ne(i,j)}(1-u_{i'}-v_{j'})\bigg)\cdot \\
 &\quad\cdot\bigg(\prod_{p<q,\,p\ne i,\,q\ne i}(u_p-u_q)\bigg)
 \bigg(\prod_{r<s,\,r\ne j,\,s\ne j}(v_r-v_s)\bigg).
\end{aligned}
\end{equation}
We will show that for every $1\le a<b\le m$, $u_a-u_b$ divides $Q$.
It will suffice to show that substituting $u_a=x$ and $u_b=x$ into 
the RHS of equation~\eqref{eq:Q} yields zero.
This substitution makes all terms of the sum zero except possibly when $i\in\{a,b\}$.
These remaining terms are equal to the common factor
\begin{align*}
x^{m-k-1}\bigg(\prod_{\substack{i'\notin\{a,b\} \\ 1\le j'\le n}}(1-u_{i'}-v_{j'})\bigg)
&\bigg(\prod_{\substack{p<q \\ p,q\notin\{a,b\}}}(u_p-u_q)\bigg)\cdot \\[1.5ex]
&\cdot\sum_{j=1}^n(-1)^jv_j^{n-\ell-1}(1-x-v_j)\prod_{j'\ne j}(1-x-v_{j'})^2
\end{align*}
times the quantity
\begin{equation}
\label{eq:pq}
(-1)^a\prod_{\substack{p<q \\ b\in\{p,q\} \\ a\notin\{p,q\}}}(u_p-u_q)
+(-1)^b\prod_{\substack{p<q \\ a\in\{p,q\} \\ b\notin\{p,q\}}}(u_p-u_q).
\end{equation}
However, since $a<b$, we see that the quantity~\eqref{eq:pq} is equal to
\begin{align*}
(-1)^a&\bigg(\prod_{p<b,\,p\ne a}(u_p-x)\bigg)\bigg(\prod_{q>b}(x-u_q)\bigg)
+(-1)^b\bigg(\prod_{p<a}(u_p-x)\bigg)\bigg(\prod_{q>a,\,q\ne b}(x-u_q)\bigg) \\[1ex]
=&\big((-1)^a(-1)^{b-a-1}+(-1)^b\big)
\bigg(\prod_{p<a}(u_p-x)\bigg)\bigg(\prod_{q>a,\,q\ne b}(x-u_q)\bigg)
=0.
\end{align*}
By symmetry, we have that $v_c-v_d$ divides the polynomial $Q$, for every $1\le c<d\le n$.
Thus
\begin{equation}
\label{eq:PQover}
P^{k,\ell}_{m,n}=\frac Q{\big(\prod_{p<q}(u_p-u_q)\big)\big(\prod_{r<s}(v_r-v_s)\big)}
\end{equation}
is a polynomial in $u_1,\dots,u_m,v_1,\dots,v_n$.
The upper bound on its degree
is easily computed from~\eqref{eq:PQover} and the expression~\eqref{eq:Q}.
\end{proof}

Proposition~\ref{prop:11v} shows $P^{k,0}_{m,1}=(1-v_1)^{m-k-1}$ and
$P^{0,\ell}_{1,n}=(1-u_1)^{n-\ell-1}$.
Here are some other examples.
\begin{examples}\rm
\label{examples:P}
\begin{align*}
P^{0,1}_{2,2} &= 1 - u_1u_2 - v_1 - 
     v_2 + v_1v_2 \\
P^{1,1}_{2,2} &= 2 - u_1 - u_2 - v_1 - 
     v_2 \\
P^{1,1}_{2,3} &= \begin{aligned}[t] &2 - 3u_1 + u_1^2 - 
     3u_2 + 4u_1u_2 - 
     u_1^2u_2 + u_2^2 - 
     u_1u_2^2 - v_1 + 
     u_1v_1 + u_2v_1 - 
     u_1u_2v_1  \\ & - v_2 + 
     u_1v_2 + u_2v_2 - 
     u_1u_2v_2 - v_3 + 
     u_1v_3 + u_2v_3 - 
     u_1u_2v_3 + 
     v_1v_2v_3 \end{aligned}\displaybreak[2] \\
P^{1,2}_{2,3} &= \begin{aligned}[t] & 3 - 3u_1 + u_1^2 - 
     3u_2 + u_1u_2 + 
     u_2^2 - 2v_1 + u_1
      v_1 + u_2v_1 - 2v_2 + 
     u_1v_2 + u_2v_2          \\ &
     + v_1v_2      
     - 2v_3 u_1v_3 + u_2v_3 + 
     v_1v_3 + v_2v_3. \end{aligned}
\end{align*}
\end{examples}

\section{Generating functions}
\label{sec:genfun}

We will now define $N(k_1,\ell_1,\dots,k_n,\ell_n)$,
but extending the definition found
in the introduction to allow $k_j,\ell_j$ to be $-1$.
The reason for this extension is to allow a recursion formula to be proved for
$N(k_1,\ell_1,\dots,k_n,\ell_n)$ that holds for all $k_j,\ell_j\in\Nats\cup\{0\}$.

\begin{defi}\rm
\label{def:N}
Let $n\in\Nats$ and $k_1,\dots,k_n,\ell_1,\dots,\ell_n\in\Nats\cup\{0,-1\}$.
If $k_1+\dots+k_n\ne\ell_1+\dots+\ell_n$ then set $N(k_1,\ell_1,\dots,k_n,\ell_n)=0$.
Suppose now $k_1+\dots+k_n=\ell_1+\dots+\ell_n$, and let $m=k_1+\dots+k_n$.
If $k_j\ge0$ and $\ell_j\ge0$ for all $j$, then let
\[
N(k_1,\ell_1,\dots,k_n,\ell_n)=(m+1)!\tau((T^*)^{k_1}T^{\ell_1}\dots(T^*)^{k_n}T^{\ell_n}).
\]
In the case $n=1, k_1=\ell_1=-1$ let $N(-1,-1)=1$.
In all other cases, (where $n\ge2$ and $k_j=-1$ or $\ell_j=-1$ for some $j$),
let $N(k_1,\ell_1,\dots,k_n,\ell_n)=0$.
\end{defi}

The following is immediate from~\cite[Proposition 9.4]{DH}.
\begin{prop}
\label{prop:Nprop}
\renewcommand{\theenumi}{\roman{enumi}}
\renewcommand{\labelenumi}{(\theenumi)}
\begin{enumerate}

\item
$N(k_1,\ell_1,\dots,k_n,\ell_n)=N(\ell_1,k_2,\ell_2,\dots,k_n,\ell_n,k_1)$.

\item
$N(k_1,\ell_1,\dots,k_n,\ell_n)=N(\ell_n,k_n,\ell_{n-1},k_{n-1},\dots,\ell_1,k_1)$.

\item
If $n\ge2$, if $k_1,\dots,k_n,\ell_1,\dots,\ell_n\ge0$ and if $k_1=0$,
then
\[
N(k_1,\ell_1,\dots,k_n,\ell_n)=N(k_2,\ell_2,\dots,k_{n-1},\ell_{n-1},k_n,\ell_n+\ell_1).
\]
\end{enumerate}
\end{prop}

\begin{thm}
\label{thm:recur}
If $n\in\Nats$ and $k_1,\dots,k_n,\ell_1,\dots,\ell_n\in\Nats\cup\{0\}$,
then
\begin{align}
N(k_1,\ell_1,&\dots,k_n,\ell_n)=\sum_{r=1}^n\sum_{1\le j(1)<\dots<j(r)\le n}
\operatorname{Nom}(\ell_1,\dots,\ell_n;j(1),\dots,j(r)) \label{eq:recur} \\
&\bigg(N(k_1,\ell_1,\dots,k_{j(1)-1},\ell_{j(1)-1},k_{j(1)}-1,
\ell_{j(r)}-1,k_{j(r)+1}\ell_{j(r)+1},\dots,k_n,\ell_n) \notag \\
&\quad\prod_{i=1}^{r-1}N(\ell_{j(i)}-1,k_{j(i)+1},\ell_{j(i)+1},\dots,
k_{j(i+1)-1},\ell_{j(i+1)-1},k_{j(i+1)}-1)\bigg), \notag
\end{align}
where $\operatorname{Nom}(\ell_1,\dots,\ell_n;j(1),\dots,j(r))$ is the multinomial coefficient
\[
\binom{\ell_1+\dots+\ell_n}
{\ell_1+\dots+\ell_{j(1)-1}+\ell_{j(r)}+\dots+\ell_n,\,
\ell_{j(1)}+\dots+\ell_{j(2)-1},\,
\dots,\,
\ell_{j(r-1)}+\dots+\ell_{j(r)-1}}.
\]
\end{thm}
\begin{proof}
If $k_1+\dots+k_n\ne\ell_1+\dots+\ell_n$ then both sides of the equality are zero.
So suppose $k_1+\dots+k_n=\ell_1+\dots+\ell_n$.
If $k_1,\dots,k_n,\ell_1,\dots,\ell_n\ge1$ then the equality
follows directly from the recursion formula~\cite[Theorem 9.5]{DH}.
To verify the equality~\eqref{eq:recur} in general, we will proceed by induction on $n$.
Let $N'(k_1,\ell_1,\dots,k_n,\ell_n)$ be the RHS of~\eqref{eq:recur}.
In the case $n=1$, we need only check that $N'(0,0)=1$, which is clear.
Assume $n\ge2$.
Writing $k_1,\ell_1,\dots,k_n,\ell_n$ on a circle, it is clear that both sides of
equation~\eqref{eq:recur} are invariant under the following permutation of order $n$:
\[
k_1\to k_2\to\dots\to k_n\to k_1,\qquad\ell_1\to\ell_2\to\dots\to\ell_n\to\ell_1.
\]
Therefore, it will be enough to consider three cases: (i) $k_1=0$, $\ell_n\ne0$,
(ii) $k_1\ne0$, $\ell_n=0$, (iii) $k_1=\ell_n=0$.

In case~(i), all terms in the sum~\eqref{eq:recur} having $j(1)=1$ are zero.
If $j(1)\ge2$, then using Proposition~\ref{prop:Nprop} we have
\begin{align*}
N(k_1,\ell_1,&\dots,k_{j(1)-1},\ell_{j(1)-1},k_{j(1)}-1,
\ell_{j(r)}-1,k_{j(r)+1}\ell_{j(r)+1},\dots,k_n,\ell_n)= \\[1ex]
&=\begin{cases}
\begin{aligned}[b]
N(&k_2,\ell_2,\dots,k_{j(1)-1},\ell_{j(1)-1},k_{j(1)}-1, \\[0.5ex]
&\ell_{j(r)}-1,k_{j(r)+1}\ell_{j(r)+1},
\dots,k_{n-1},\ell_{n-1},k_n,\ell_n+\ell_1)\end{aligned}
&\text{if }j(r)<n \\[1ex]
N(k_2,\ell_2,\dots,k_{j(1)-1},\ell_{j(1)-1},k_{j(1)}-1,
\ell_{j(r)}+\ell_1-1)
&\text{if }j(r)=n.
\end{cases}
\end{align*}
Therefore,
\begin{align*}
N'(k_1,\ell_1,\dots,k_n,\ell_n)&=N'(k_2,\ell_2,\dots,k_{n-1},\ell_{n-1},k_n,\ell_n+\ell_1) \\
&=N(k_2,\ell_2,\dots,k_{n-1},\ell_{n-1},k_n,\ell_n+\ell_1) \\
&=N(0,\ell_1,k_2,\ell_2,\dots,k_n,\ell_n),
\end{align*}
where the second equality is by the induction hypothesis
and the third is from Proposition~\ref{prop:Nprop}.

In case~(ii), ($k_1\ne0$, $\ell_n=0$), one similarly shows
\begin{align*}
N'(k_1,\ell_1,\dots,k_n,\ell_n)&=N'(k_1+k_n,\ell_1,k_2,\ell_2,\dots,k_{n-1},\ell_{n-1}) \\
&=N(k_1+k_n,\ell_1,k_2,\ell_2,\dots,k_{n-1},\ell_{n-1})\\
&=N(k_1,\ell_1,,\dots k_{n-1},\dots,k_{n-1},\ell_{n-1},k_n,0).
\end{align*}

In case~(iii), ($k_1=\ell_n=0$), given $1\le r\le n$ and $1<j(1)<\dots<j(r)\le n$, let
$t(j(1),\dots,j(r))$ be the corresponding term in the sum~\eqref{eq:recur} for
$N'(k_1,\ell_1,\dots,k_n,\ell_n)$, and, if $j(r)<n$,
let $\ttild(j(1),\dots,j(r))$ be the term corresponding to
$1\le j(1)<\dots<j(r)\le n-1$
in the analogous sum for $N'(k_2,\ell_2,\dots,k_{n-1},\ell_{n-1},k_n,\ell_1)$.
If $j(1)>1$ and $j(r)<n$, then using Proposition~\ref{prop:Nprop} we have
\begin{align*}
&N(k_1,\ell_1,\dots,k_{j(1)-1},\ell_{j(1)-1},k_{j(1)}-1,\ell_{j(r)}-1,
k_{j(r)+1},\ell_{j(r)+1},\dots,k_n,\ell_n)= \\[1ex]
&\,N(k_2,\ell_2,\dots,k_{j(1)-1},\ell_{j(1)-1},k_{j(r)}-1,\ell_{j(r)}-1,
k_{j(r)+1},\ell_{j(r)+1},\dots,k_{n-1},\ell_{n-1},k_n,\ell_1),
\end{align*}
which, taken in the sum~\eqref{eq:recur}, shows
\[
t(j(1),\dots,j(r))
=\ttild(j(1)-1,\dots,j(r)-1).
\]
If $j(1)=1$ and $j(r)<n$ or if $j(1)>1$ and $j(r)=n$, then
$t(j(1),\dots,j(r))=0$.
If $j(1)=1$ and $j(r)=n$ then since
\begin{align*}
N(-1,-1)N(\ell_1-1,&k_2,\ell_2,\dots,k_{j(2)-1},\ell_{j(2)-1},k_{j(2)}-1) \\
&=N(k_2,\ell_2,\dots,k_{j(2)-1},\ell_{j(2)-1},k_{j(2)}-1,\ell_1-1),
\end{align*}
we have
\[
t(j(1),\dots,j(r)=\ttild(j(2)-1,\dots,j(r-1)-1,n-1).
\]
Therefore,
\begin{align*}
N'(k_1,\ell_1,\dots,k_n,\ell_n)&=N'(k_2,\ell_2,\dots,k_{n-1},\ell_{n-1},k_n,\ell_1) \\
&=N(k_2,\ell_2,\dots,k_{n-1},\ell_{n-1},k_n,\ell_1) \\
&=N(0,\ell_1,k_2,\ell_2,\dots,k_{n-1},\ell_{n-1},k_n,0).
\end{align*}
\end{proof}

Recall the generating functions $F_n$, defined in the introduction.
The following theorem is a recursion formula for these, using the operation
$\odot=\odot_2$ introduced in~\S\ref{sec:fps}.
Here, as promised, we are using objects defined in~\S\ref{sec:fps} with $N=2$.

\begin{thm}
\label{thm:Frecur}
If $n\in\Nats$ and $n\ge2$, then
\begin{align}
&(1-z_1w_1-\dots-z_nw_n)F_n(z_1,w_1,\dots,z_n,w_n)=
\sum_{r=2}^n\sum_{1\le j(1)<\dots<j(r)\le n} \label{eq:Frecur} \\
&\big(z_{j(1)}w_{j(r)}\Ft_{j(1)+n-j(r)}(z_1,w_1,\dots,z_{j(1)-1},w_{j(1)-1},z_{j(1)},
w_{j(r)},z_{j(r)+1},w_{j(r)+1},\dots,z_n,w_n)\big) \notag \\[1ex]
&\odot\bigg(\operatornamewithlimits{\odot}_{i=1}^{r-1}\big(w_{j(i)}z_{j(i+1)}
\Ft_{j(i+1)-j(i)}(w_{j(i)},z_{j(i)+1},w_{j(i)+1},\dots,z_{j(i+1)-1},w_{j(i+1)-1},z_{j(i+1)})\big)\bigg),
\notag
\end{align}
where
\[
\Ft_j(x_1,y_1,\dots,x_j,y_j):=\begin{cases}
F_j(x_1,y_1,\dots,x_j,y_j)&\text{if }j\ge2 \\
F_1(x_1,y_1)/(x_1y_1)&\text{if }j=1.
\end{cases}
\]
\end{thm}
\begin{proof}
We will use the recursion formula~\eqref{eq:recur} for the coefficients
in the definition~\eqref{eq:Fdef} of $F_n$.
Since $n\ge2$, if $j\in\{1,\dots,n\}$, then
\begin{align*}
&\sum_{\substack{k_1,\dots,k_n\ge0 \\ \ell_1,\dots,\ell_n\ge0}}
N(k_1,\ell_1,\dots,k_{j-1},\ell_{j-1},k_j-1,\ell_j-1,k_{j+1},\ell_{j+1},\dots,k_n,\ell_n)
z_1^{k_1}w_1^{\ell_1}\dots z_n^{k_n}w_n^{\ell_n} \\[1ex]
&\qquad\qquad=\quad z_jw_jF_n(z_1,w_1,\dots,z_n,w_n).
\end{align*}
This takes care of all terms in the sum~\eqref{eq:recur} with $r=1$,
and these together with the LHS of equation~\eqref{eq:recur}
give the LHS of formula~\eqref{eq:Frecur}.
Let $2\le r\le n$ and $1\le j(1)<\dots<j(r)\le n$.
Then
\begin{align}
&\sum_{\substack{k_1,\ldots,k_n\ge0 \\ \ell_1,\ldots,\ell_n\ge0}}
\operatorname{Nom}(\ell_1,\ldots,\ell_n;j(1),\ldots,j(r))\cdot \notag \\
&\qquad\quad\cdot\bigg(N(k_1,\ell_1,\ldots,k_{j(1)-1},\ell_{j(1)-1},k_{j(1)}-1,
\ell_{j(r)}-1,k_{j(r)+1}\ell_{j(r)+1},\ldots,k_n,\ell_n) \notag \\
&\qquad\qquad\quad\prod_{i=1}^{r-1}N(\ell_{j(i)}-1,k_{j(i)+1},\ell_{j(i)+1},\ldots,
k_{j(i+1)-1},\ell_{j(i+1)-1},k_{j(i+1)}-1)\bigg)\cdot \notag \\[1ex]
&\qquad\qquad\cdot z_1^{k_1}w_1^{\ell_1}\ldots z_n^{k_n}w_n^{\ell_n}
 \notag \displaybreak[2] \\[2ex]
=&\sum_{p_1,\dots,p_r\ge0}\binom{p_1+\dots+p_r}{p_1,\,\dots,\,p_r} \notag \\
&\quad\Bigg(\quad\bigg(\sum_{
 \substack{k_1+\ell_1+\dots+k_{j(1)-1}+\ell_{j(1)-1}+k_{j(1)}+ \\
           +\ell_{j(r)}+k_{j(r)+1}+\ell_{j(r)+1}+\dots+k_n+\ell_n=p_1} } \label{eq:sum1} \\[1ex]
&\qquad\qquad
N(k_1,\ell_1,\dots,k_{j(1)-1},\ell_{j(1)-1},k_{j(1)}-1,\ell_{j(r)}-1,
k_{j(r)+1},\ell_{j(r)+1},\dots,k_n,\ell_n)\cdot \notag \\[1ex]
&\qquad\qquad\qquad\qquad
\cdot z_1^{k_1}w_1^{\ell_1}\dots z_{j(1)-1}^{k_{j(1)-1}}w_{j(1)-1}^{\ell_{j(1)-1}}
z_{j(1)}^{k_{j(1)}}w_{j(r)}^{\ell_{j(r)}}z_{j(r)+1}^{k_{j(r)+1}}w_{j(r)+1}^{\ell_{j(r)+1}}
\dots z_n^{k_n}w_n^{\ell_n}\bigg)\cdot \notag \displaybreak[2] \\[1ex]
&\quad\cdot\bigg(\prod_{i=1}^{r-1}\sum_{\substack{
\ell_{j(i)}+k_{j(i)+1}+\ell_{k(i)+1}+\dots+ \\+k_{j(i+1)-1}+\ell_{j(i+1)-1}+k_{j(i+1)}=p_{i+1} } }
\label{eq:sumi} \\[1ex]
&\qquad\qquad
N(\ell_{j(i)}-1,k_{j(i)+1},\ell_{j(i)+1},\dots,k_{j(i+1)-1},\ell_{j(i+1)-1},k_{j(i+1)}-1)\cdot \notag \\
&\qquad\qquad\qquad\qquad\qquad
\cdot w_{j(i)}^{\ell_{j(i)}}z_{j(i)+1}^{k_{j(i)+1}}w_{j(i)+1}^{\ell_{j(i)+1}}\dots
z_{j(i+1)-1}^{k_{j(i+1)-1}}w_{j(i+1)-1}^{\ell_{j(i+1)-1}}z_{j(i+1)}^{k_{j(i+1)}}\bigg)\quad\Bigg) \notag
\end{align}
Consider a summation beginning at~\eqref{eq:sumi} for some fixed $i$.
If $j(i+1)-j(i)\ge2$ then all terms with either $\ell_{j(i)}=0$ or $k_{j(i+1)}=0$ vanish
and the summation is equal to the $2p_{i+1}$--homogeneous part of
\[
w_{j(i)}z_{j(i+1)}F_{j(i+1)-j(i)}(w_{j(i)},z_{j(i)+1},w_{j(i)+1},
\dots,z_{j(i+1)-1},w_{j(i+1)-1},z_{j(i+1)}).
\]
On the other hand, if $j(i+1)-j(i)=1$, then, since $N(k-1,\ell-1)=N(k,\ell)$ for all
$k,\ell\in\Nats\cup\{0\}$, this summation is equal to
the $2p_{i+1}$--homogenous part of $F_1(w_{j(i)},z_{j(i+1)})$.
In either case, the summation is equal to the $2p_{i+1}$--homogeneous part of
\[
w_{j(i)}z_{j(i+1)}\Ft_{j(i+1)-j(i)}(w_{j(i)},z_{j(i)+1},w_{j(i)+1},
\dots,z_{j(i+1)-1},w_{j(i+1)-1},z_{j(i+1)}).
\]
Similarly, we see that the summation beginning at~\eqref{eq:sum1} is equal to
 the $2p_1$--homogeneous part of
\[
z_{j(1)}w_{j(r)}\Ft_{j(1)+n-j(r)}(z_1,w_1,\dots,z_{j(1)-1},w_{j(1)-1},z_{j(1)},
w_{j(r)},z_{j(r)+1},w_{j(r)+1},\dots,z_n,w_n).
\]
Now using the definition of the operation $\odot$
and summing over $2\le r\le n$ and all $1\le j(1)<\cdots<j(r)\le n$,
one proves the recursion formula~\eqref{eq:recur}.
\end{proof}

It is straightforward from the recursion formula that
\[
F_1(z_1,w_1)=\frac1{1-z_1w_1}.
\]
Starting from this, one can (in principle) compute $F_n$ for arbitrary given $n$,
using the recursion formula~\eqref{eq:recur}.
Below are the results for $n=2,3,4$.
We were motivated to find and write down $F_4$, although
it is rather long, because, while we were able to find $F_2$ and $F_3$ using
ad hoc methods to work out, in the notation of Proposition~\ref{prop:11v},
\[
\frac1{1-u_1}\odot\frac1{1-v}\qquad\text{and}\qquad
\frac a{(1-u_1)(1-u_2)}\odot\frac1{1-v},
\]
we were unable to find, in the notation of Proposition~\ref{prop:adtb},
\[
\frac a{(1-u_1)(1-u_2)}\odot\frac b{(1-v_1)(1-v_2)},
\]
which is needed in computing $F_4$,
until we discovered the $\EEu$--transform.

\begin{examples}\rm
($n=2$):
\begin{align*}
F_2(z_1,w_1,z_2,w_2)&=\frac1{1-u_1^{(2)}}\big(F_1(z_1,w_2)\odot F_1(w_1,z_2)\big) \\
&=\frac1{1-u_1^{(2)}}\bigg(\frac1{1-z_1w_2}\odot \frac1{1-w_1z_2}\bigg) \\
&=\frac1{(1-u_1^{(2)})(1-u_2^{(2)})}
\end{align*}
where $u_1^{(2)}=z_1w_1+z_2w_2$ and $u_2^{(2)}=z_1w_2+z_2w_1$.

\vskip2ex
\noindent
($n=3$):
\begin{align*}
&F_3(z_1,w_1,z_2,w_2,z_3,w_3)= \\
&=\frac1{1-u_1^{(3)}}\begin{aligned}[t]
\bigg(&
\big(z_1w_2F_2(z_1,w_2,z_3,w_3)\big)\odot F_1(w_1,z_2) \\
&\qquad+\;F_1(z_1,w_3)\odot\big(w_1z_3F_2(w_1,z_2,w_2,z_3)\big) \\[1.5ex]
&\qquad\qquad\qquad+\;\big(z_2w_3F_2(z_1,w_1,z_2,w_3)\big)\odot F_1(w_2,z_3) \\
&\qquad\qquad\qquad\qquad\qquad+\;F_1(z_1,w_3)\odot F_1(w_1,z_2)\odot F_1(w_2,z_3)\bigg)
\end{aligned} \\
&=\frac1{(1-u_1^{(3)})(1-u_2^{(3)})}\bigg(1+
\frac{z_1w_2}{1-u_3^{(3)}}+\frac{z_2w_3}{1-u_4^{(3)}}
+\frac{z_3w_1}{1-u_5^{(3)}}\bigg)
\end{align*}
where
\begin{alignat*}{2}
u_1^{(3)}&=z_1w_1+z_2w_2+z_3w_3 \qquad & u_2^{(3)}&=z_1w_3+z_2w_1+z_3w_2 \\
u_3^{(3)}&=z_1w_2+z_2w_1+z_3w_3 \qquad & u_4^{(3)}&=z_1w_1+z_2w_3+z_3w_2 \\
u_5^{(3)}&=z_1w_3+z_2w_2+z_3w_1.
\end{alignat*}

\vskip2ex
\noindent
($n=4$):
\begin{align*}
F_4(z_1,&w_1,z_2,w_2,z_3,w_3,z_4,w_4)= \\
&=\frac1{1-u^{(4)}_1}\bigg(
( z_1w_2 F_3(z_1, w_2, z_3, w_3, z_4, w_4))\odot F_1(w_1, z_2) \\
&\qquad\qquad\qquad+\quad( z_1 w_3 F_2(z_1, w_3, z_4, w_4))\odot ( w_1z_3 F_2(w_1, z_2, w_2, z_3)) \\[1.5ex]
&\qquad\qquad\qquad+\quad  F_1(z_1, w_4)\odot (  w_1 z_4 F_3(w_1, z_2, w_2, z_3, w_3, z_4))
\displaybreak[2] \\[1.5ex]
&\qquad\qquad\qquad+\quad (  z_2w_3  F_3(z_1, w_1, z_2, w_3, z_4, w_4))\odot F_1(w_2, z_3)
\displaybreak[2] \\[1.5ex]
&\qquad\qquad\qquad+\quad ( z_2 w_4 F_2(z_1, w_1, z_2, w_4) )\odot ( w_2z_4 F_2(w_2, z_3, w_3, z_4))
\displaybreak[2] \\[1.5ex]
&\qquad\qquad\qquad+\quad (  z_3 w_4 F_3(z_1, w_1, z_2, w_2, z_3, w_4))\odot F_1(w_3, z_4)
\displaybreak[2] \\[1.5ex]
&\qquad\qquad\qquad+\quad ( z_1 w_3 F_2(z_1, w_3, z_4, w_4) )\odot F_1(w_1, z_2) \odot F_1(w_2, z_3)
\displaybreak[2] \\[1.5ex]
&\qquad\qquad\qquad+\quad F_1(z_1, w_4) \odot F_1(w_1, z_2)\odot ( w_2 z_4 F_2(w_2, z_3, w_3, z_4))
\displaybreak[2] \\[1.5ex]
&\qquad\qquad\qquad+\quad F_1(z_1, w_4)\odot ( w_1 z_3  F_2(w_1, z_2, w_2, z_3)) \odot  F_1(w_3, z_4)
\displaybreak[2] \\[1.5ex]
&\qquad\qquad\qquad+\quad ( z_2 w_4 F_2(z_1, w_1, z_2, w_4))\odot F_1(w_2, z_3)\odot F_1(w_3, z_4)
 \\
&\qquad\qquad\qquad+\quad F_1(z_1, w_4)) \odot  F_1(w_1, z_2) \odot F_1(w_2, z_3) \odot F_1(w_3, z_4)\bigg)
\displaybreak[3] \\
&=\frac1{(1 - u^{(4)}_1)(1 - u^{(4)}_2)}\bigg(1+\frac{z_3w_4}{1 - u^{(4)}_3}
+ \frac{z_2w_3}{1 - u^{(4)}_5} 
+ \frac{z_4w_1}{1 - u^{(4)}_6}
+ \frac{z_1w_2}{1 - u^{(4)}_8}+ \\[1ex]
&\qquad+ \frac{z_2w_4}{1 - u^{(4)}_9}
+ \frac{z_3w_1}{1 - u^{(4)}_{10}}
+ \frac{z_4w_2}{1 - u^{(4)}_{11}}
+ \frac{z_1w_3}{1 - u^{(4)}_{12}}
+ \frac{z_1z_3w_2w_4}{(1 - u^{(4)}_3)(1 - u^{(4)}_{14})}+ \displaybreak[2] \\[1ex]
&\qquad+ \frac{z_2z_3w_4^2}{(1 - u^{(4)}_3)(1 - u^{(4)}_9)}
+ \frac{z_3^2w_1w_4}{(1 - u^{(4)}_3)(1 - u^{(4)}_{10})}
+ \frac{z_2z_4w_2w_4}{(1 - u^{(4)}_4)(1 - u^{(4)}_9)} + \displaybreak[2] \\[1ex]
&\qquad+ \frac{z_2z_4w_2w_4}{(1 - u^{(4)}_4)(1 - u^{(4)}_{11})}
+ \frac{z_2^2w_3w_4}{(1 - u^{(4)}_5)(1 - u^{(4)}_9)}
+ \frac{z_1z_2w_3^2}{(1 - u^{(4)}_5)(1 - u^{(4)}_{12})} +\displaybreak[2]\\[1ex]
&\qquad+ \frac{z_2z_4w_1w_3}{(1 - u^{(4)}_5)(1 - u^{(4)}_{13})}
+ \frac{z_3z_4w_1^2}{(1 - u^{(4)}_6)(1 - u^{(4)}_{10})}
+ \frac{z_4^2w_1w_2}{(1 - u^{(4)}_6)(1 - u^{(4)}_{11})} + \displaybreak[2]\\[1ex]
&\qquad+ \frac{z_2z_4w_1w_3}{(1 - u^{(4)}_6)(1 - u^{(4)}_{13})}
+ \frac{z_1z_3w_1w_3}{(1 - u^{(4)}_7)(1 - u^{(4)}_{10})}
+ \frac{z_1z_3w_1w_3}{(1 - u^{(4)}_7)(1 - u^{(4)}_{12})}+ \displaybreak[2]\\[1ex]
&\qquad+ \frac{z_1z_4w_2^2}{(1 - u^{(4)}_8)(1 - u^{(4)}_{11})} 
+ \frac{z_1^2w_2w_3}{(1 - u^{(4)}_8)(1 - u^{(4)}_{12})}
+ \frac{z_1z_3w_2w_4}{(1 - u^{(4)}_8)(1 - u^{(4)}_{14})}\bigg),
\end{align*}
where
\begin{alignat*}{2}
u^{(4)}_1 &= z_1 w_1 + z_2 w_2 + z_3 w_3 + z_4 w_4 \qquad&
u^{(4)}_2 &= z_1 w_4 + z_2 w_1 + z_3 w_2 + z_4 w_3 \\
u^{(4)}_3 &= z_1 w_1 + z_2 w_2 + z_3 w_4 + z_4 w_3 &
u^{(4)}_4 &= z_1 w_1 + z_2 w_4 + z_3 w_3 + z_4 w_2 \\
u^{(4)}_5 &= z_1 w_1 + z_2 w_3 + z_3 w_2 + z_4 w_4 &
u^{(4)}_6 &= z_1 w_4 + z_2 w_2 + z_3 w_3 + z_4 w_1 \\
u^{(4)}_7 &= z_1 w_3 + z_2 w_2 + z_3 w_1 + z_4 w_4 &
u^{(4)}_8 &= z_1 w_2 + z_2 w_1 + z_3 w_3 + z_4 w_4 \\
u^{(4)}_9 &= z_1 w_1 + z_2 w_4 + z_3 w_2 + z_4 w_3 &
u^{(4)}_{10} &= z_1 w_4 + z_2 w_2 + z_3 w_1 + z_4 w_3 \\
u^{(4)}_{11} &= z_1 w_4 + z_2 w_1 + z_3 w_3 + z_4 w_2 &
u^{(4)}_{12} &= z_1 w_3 + z_2 w_1 + z_3 w_2 + z_4 w_4 \\
u^{(4)}_{13} &= z_1 w_4 + z_2 w_3 + z_3 w_2 + z_4 w_1 &
u^{(4)}_{14} &= z_1 w_2 + z_2 w_1 + z_3 w_4 + z_4 w_3.
\end{alignat*}
\end{examples}

\begin{thm}
\label{thm:rat}
For every $n\in\Nats$ with $n\ge2$,
$F_n(z_1,w_1,\dots,z_n,w_n)$ is a rational function in variables
$z_1,\dots,z_n,w_1,\dots,w_n$.
It is, moreover, a sum of terms of the form
\[
\frac{z_{i_1}\dots z_{i_k}\,w_{j_1}\dots w_{j_k}}
{(1-u_1)\dots(1-u_\ell)}
\]
for some $\ell\in\{2,\dots,n\}$, $k\in\{0,1,\dots,\ell-1\}$,
$i_1,\dots,i_k,j_1,\dots,j_k\in\{1,\dots,n\}$ and for $u_1,\dots,u_\ell$
of the form
\[
u_p=\sum_{i=1}^nz_iw_{\sigma_p(i)}
\]
for distinct permutations $\sigma_1,\dots,\sigma_\ell$ of $\{1,\dots,n\}$.
\end{thm}
\begin{proof}
This follows by induction on $n$ from the recursion formula~\eqref{eq:recur}
and Proposition~\ref{prop:adtb}.
\end{proof}

Examining the power series of $F_n$, for each $n$ the conjecture~\eqref{eq:conj}
yields a conjectured identity involving multinomial coefficients.
For example, when $n=3$ we get the conjecture
\begin{align*}
3^{3p}=&
\sum_{\substack{j,k\ge0\\j+k=p}}\binom{3j}{j,j,j}\binom{3k}{k,k,k}+ \\[1ex]
&+3\sum_{\substack{j,k,\ell\ge0\\j+k+\ell=p-1}}\;
\sum_{\substack{k',\ell'\ge0\\k'+\ell'=k+\ell+1}}\;
\sum_{\substack{j',\ell''\ge0\\j'+\ell''=j+\ell+1}}
\binom{2j+j'}{j,\,j,\,j'}\binom{2k+k'}{k,\,k,\,k'}\binom{\ell+\ell'+\ell''}{\ell,\,\ell',\,\ell''}
\end{align*}
for all $p\in\Nats$.

\section{Contour integration}
\label{sec:ci}

It follows from Theorem~\ref{thm:rat} that $F_n(z_1,w_1,\ldots,z_n,w_n)$
is a holomorphic function of complex variables $z_1,w_1,\ldots,z_n,w_n$
in a suitably small ball around the origin.
It is possible, at least in principle,
to verify the conjecture mentioned at equation~\eqref{eq:conj}
in the introduction for a given value of $n$ from the generating function $F_n$
by performing $2n-1$ contour integrations.
For example, letting
\[
G_n(x_1,\dots,x_n)=\sum_{k_1,\dots,k_n\ge0}N(k_1,k_1,k_2,k_2,\dots,k_n,k_n)x_1^{k_1}\dots x_n^{k_n},
\]
for suitably small $\eps>0$ we have
\begin{multline*}
G_n(y_1^2,\ldots,y_n^2)= \\
=\frac1{(2\pi i)^n}\int_{|\zeta_1|=\eps}\dots\int_{|\zeta_n|=\eps}
\zeta_1^{-1}\dots\zeta_n^{-1}
F_n(\zeta_1y_1,\zeta_1^{-1}y_1,\ldots,\zeta_ny_n,\zeta_n^{-1}y_n)\dif\zeta_n\dots\dif\zeta_1.
\end{multline*}
Letting
\[
H_n(a)=\sum_{k=0}^\infty N(\underset{2n\text{ times}}{\underbrace{k,k,\ldots,k,k}})a^k,
\]
we have
\begin{align*}
&H_n(b^n)= \\
&\;=\frac1{(2\pi i)^{n-1}}\int_{|\eta_2|=\eps}\dots\int_{|\eta_n|=\eps}
\eta_2^{-1}\dots\eta_n^{-1}
G_n(\eta_2\dots\eta_n b,\eta_2^{-1}b,\eta_3^{-1}b,\dots,\eta_n^{-1}b)\dif\eta_n\dots\dif\eta_2.
\end{align*}
The conjecture~\eqref{eq:conj} is that
\[
H_n(a)=\frac1{1-n^na}
\]
for every $n\in\Nats$.

Although the counjecture can be verified without difficulty
in the case $n=2$ using combinatorial methods,
we will illustrate the method described above in checking it.
However, even for $n=3$, we have been unable to prove the conjecture using
the method of contour integration.

We have
\begin{align*}
G_2^{(1)}(y_1,z_2,w_2):=&\frac1{2\pi i}
\int_{|\zeta|=\eps}F_2(\zeta y_1,\zeta^{-1}y_1,z_2,w_2)\zeta^{-1}\dif\zeta \\
=&\frac1{2\pi i(1-y_1^2-z_2w_2)}\int_{|\zeta|=\eps}\frac1{\zeta-\zeta^2y_1w_2-y_yz_2}\dif\zeta
\end{align*}
and
$\zeta-\zeta^2y_1w_2-y_1z_2=y_1w_1(r_1-\zeta)(\zeta-r_2)$
with
\[
r_1,r_2=\frac{1\pm\sqrt{1-4y_1^2z_2w_2}}{2y_1w_2}.
\]
If $|y_1|,\,|z_2|,\,|w_2|\,<\delta$ for $\delta>0$ small enough, $r_1$ can be forced arbitrarily close
to zero and $r_2$ can be forced close to $\infty$.
By the residue theorem,
\begin{align*}
G_2^{(1)}(y_1,z_2,w_2)&=\frac1{(1-y_1^2-z_2w_2)y_1w_2(r_1-r_2)} \\[1ex]
&=\frac1{(1-y_1^2-z_2w_2)\sqrt{1-4y_1^2z_2w_2}}.
\end{align*}
Therefore,
\[
G_2(x_1,x_2)=\frac1{(1-x_1-x_2)\sqrt{1-4x_1x_2}}.
\]
Finally,
\begin{align*}
H_2(b^2)&=\frac1{2\pi i}\int_{|\eta|=1}\eta^{-1}G_2(\eta b,\eta^{-1}b)\dif\eta \\
&=\frac1{2\pi i\sqrt{1-4b^2}}\int_{|\eta|=1}\frac1{\eta-b-\eta^2b}\dif\eta \\
&=\frac1{1-4b^2}.
\end{align*}

\bibliographystyle{plain}

\end{document}